\documentclass{article}

\usepackage{amsmath,amsfonts,amssymb,stmaryrd,enumerate}
\usepackage{fullpage}

\newtheorem{innercustomthm}{Theorem}
\newenvironment{thm}[1]
  {\renewcommand\theinnercustomthm{#1}\innercustomthm}
  {\endinnercustomthm}
  
\providecommand{\charf}[1]{\llbracket #1 \rrbracket}

\title{A comment on \\\emph{Intersecting Families of Permutations}}
\author{Yuval Filmus}

\begin{document}

\maketitle
\begin{abstract}
Ellis, Friedgut and Pilpel~\cite{EFP} prove that for large enough $n$, a $t$-intersecting family of permutations contains at most $(n-t)!$ permutations. Their main theorem also states that equality holds only for $t$-cosets. We show that their proof of the characterization of extremal families is wrong. However, the characterization follows from a paper of Ellis~\cite{Ellis-stability}, as mentioned already by Ellis, Friedgut and Pilpel.
\end{abstract}

\section{Introduction} \label{sec:introduction}

The classical Erd\H{o}s--Ko--Rado theorem states that when $n > 2k$, an intersecting family of $k$-subsets of $[n]:=\{1,\ldots,n\}$ contains at most $\binom{n-1}{k-1}$ sets, and moreover this is achieved only by \emph{stars}, consisting of all sets containing a specific point. Wilson~\cite{Wilson} extended this result to $t$-intersecting families: he showed that when $n > (t+1)(k-t+1)$, a $t$-intersecting family of $k$-subsets of $[n]$ contains at most $\binom{n-t}{k-t}$ sets, and moreover this is achieved only by \emph{$t$-stars}, consisting of all sets containing $t$~specific points.

The Erd\H{o}s--Ko--Rado theorem has been extended in many different directions, one of them to other domains. One of the most intriguing domains is that of \emph{permutations}. We say that two permutations $\pi,\sigma \in S_n$ are $t$-intersecting if there exist $t$ points $i_1 < \cdots < i_t$ such that $\pi(i_1) = \sigma(i_1),\ldots,\pi(i_t) = \sigma(i_t)$. A $t$-intersecting family of permutations is one in which every two permutations $t$-intersect.

A simple partitioning argument shows that any $1$-intersecting family of permutations contains at most $(n-1)!$ permutations. This bound is only achieved by \emph{cosets} $T_{ij} = \{ \pi \in S_n : \pi(i) = j \}$, though this is surprisingly hard to show (different proofs appear in~\cite{CameronKu,LaroseMalvenuto,GodsilMeagher,EFP}); see also the corresponding stability result of Ellis~\cite{Ellis-cameron-ku}.

In a groundbreaking paper, Ellis, Friedgut and Pilpel~\cite{EFP} showed that for any $t$ there exists $C_t$ such that for all $n \geq C_t$, any $t$-intersecting family of permutations contains at most $(n-t)!$ permutations, and furthermore this is only achieved by \emph{$t$-cosets} $T_{\alpha_1\ldots \alpha_t \mapsto \beta_1\ldots \beta_t} = \{ \pi \in S_n : \pi(\alpha_1) = \beta_1,\ldots,\pi(\alpha_t) = \beta_t \}$. 

Unfortunately, the proof in~\cite{EFP} that $t$-cosets are the only $t$-intersecting families of permutations of maximum size is wrong for $t > 1$, as we indicate below. Fortunately, the main result of Ellis~\cite{Ellis-stability} implies that this theorem is correct; this is mentioned in~\cite{EFP} as an alternative proof of the characterization of $t$-intersecting families of maximum size. The same problem affects the characterization of the optimal families for setwise-$t$-intersecting families of permtuations in Ellis~\cite{Ellis-setwise}, but according to Ellis (private communication), the result can be recovered using methods similar to those in~\cite{Ellis-stability}.

\smallskip

The remainder of this note is structured as follows. In Section~\ref{sec:efp} we outline the relevant portion of the arguments of~\cite{EFP}. In Section~\ref{sec:thm27} and Section~\ref{sec:thm29} we present counterexamples to two results in~\cite{EFP}. We close the note with some final remarks in Section~\ref{sec:remarks}.

\paragraph{Acknowledgements} I thank David Ellis for several useful discussions.

\pagebreak

\section{The argument of Ellis--Friedgut--Pilpel} \label{sec:efp}

Here is the main theorem of~\cite{EFP}.

\begin{thm}{3}
 For any $k \in \mathbb{N}$, and any $n$ sufficiently large depending on $k$, if $I \subset S_n$ is $k$-intersecting, then $|I| \leq (n-k)!$. Equality holds if and only if $I$ is a $k$-coset of $S_n$.
\end{thm}

This first step toward proving Theorem~3 is the following result, which combines parts of Theorem~5 and Theorem~6 of~\cite{EFP}.

\begin{thm}{5/6}
 For any $k \in \mathbb{N}$ and any $n$ sufficiently large depending on $k$, if $I \subset S_n$ is $k$-intersecting, then $|I| \leq (n-k)!$. Moreover, if $I \subset S_n$ is a $k$-intersecting family of size $(n-k)!$, then $1_I \in V_k$.
\end{thm}

Here $1_I$ is the characteristic function of the family $I$, and $V_k$ is the linear span of the characteristic functions of the $k$-cosets. Given Theorem~5/6, to complete the proof of Theorem~3 we need to show that if $I$ is a $k$-intersecting family of size $(n-k)!$, then $I$ is a $k$-coset. The authors deduce this from the following result (called Theorem~8 in the introduction of~\cite{EFP}), in which $\mathcal{A}_k$ consists of all ordered $k$-tuples of distinct numbers in $[n]$.

\begin{thm}{27}
 Let $f \in V_k$ be nonnegative. Then there exist nonnegative coefficients $(b_{\alpha,\beta})_{\alpha,\beta \in \mathcal{A}_k}$ such that $f = \sum b_{\alpha,\beta} 1_{T_{\alpha \mapsto \beta}}$. Furthermore, if $f$ is Boolean, then $f$ is the characteristic function of a disjoint union of $k$-cosets.
\end{thm}

Unfortunately, this theorem is wrong, as we show in Section~\ref{sec:thm27}. However, the special case $k=1$, proved in~\cite{EFP} separately as Theorem~28, is correct.

The proof of Theorem~27 relies on two theorems, Theorem~29 and Theorem~30. Of these, Theorem~29 is incorrect, as we show in Section~\ref{sec:thm29}, in which we also indicate the mistake in the proof.

\section{Counterexample to Theorem 27} \label{sec:thm27}

From now on we denote the characteristic function of a set $S$ by $\charf{S}$ instead of $1_S$, to increase legibility.

Let $n \geq 6$, and consider the following function:
\[
 f(\pi) = \charf{\pi(\{1,2,3\}) = \{1,2,3\} \text{ or } \pi(\{1,2,3\}) \cap \{1,2,3\} = \emptyset}.
\]
This is clearly a Boolean function. We claim that it belongs to $V_2$. Indeed, we can rewrite it in the following form, which makes it clear that it belongs to $V_2$:
\begin{gather*}
 f(\pi) = 1 - \charf{\pi(1) \in \{1,2,3\}} - \charf{\pi(2) \in \{1,2,3\}} - \charf{\pi(3) \in \{1,2,3\}} \\ + \charf{\pi(\{1,2\}) \subset \{1,2,3\}} + \charf{\pi(\{1,3\}) \subset \{1,2,3\}} + \charf{\pi(\{2,3\}) \subset \{1,2,3\}}.
\end{gather*}
If Theorem~27 were true, then the support of $f$ would be the disjoint union of $2$-cosets. In particular, it would contain some $2$-coset. However, it cannot contain any $2$-coset, since knowing the image of any two points in a permutation $\pi$ is not enough to certify that $f(\pi) = 1$. (We leave it to the reader to perform the necessary case analysis.)

This counterexample can be generalized to any even $k \geq 2$ and $n \geq 2(k+1)$:
\begin{gather*}
 \charf{\pi(\{1,\ldots,k+1\}) = \{1,\ldots,k+1\} \text{ or } \pi(\{1,\ldots,k+1\}) \cap \{1,\ldots,k+1\} = \emptyset} = \\
 \sum_{S \subsetneq \{1,\ldots,k+1\}} (-1)^{|S|} \charf{\pi(S) \subset \{1,\ldots,k+1\}}.
\end{gather*}
This is a function in $V_k$ whose support contains no $k$-coset.

We describe other counterexamples in Section~\ref{sec:remarks}.

\section{Counterexample to Theorem 29} \label{sec:thm29}

The proof of Theorem~27 relies on a claimed generalization of Birkhoff's theorem attributed to Benabbas, Friedgut, and Pilpel. Birkhoff's theorem states that every bistochastic matrix is a convex combination of permutation matrices. Theorem~29 of~\cite{EFP} purports to generalize Birkhoff's theorem to $k$~dimensions. Since our counterexample is two-dimensional, we will concentrate on the case $k=2$.

A 2-bistochastic matrix $M$ is an $n(n-1) \times n(n-1)$ matrix whose rows and columns are indexed by pairs of indices, such that the following two properties hold:
\begin{enumerate}[(S1)]
 \item There exist an $n\times n$ bistochastic matrix $R = (r_{i,j})$ and $n^2$ bistochastic matrices $M_{i,j}$ of dimension $(n-1)\times(n-1)$ whose rows and vertices are indexed by $[n] \setminus \{i\}, [n] \setminus \{j\}$ (respectively) such that $M((i,i'),(j,j')) = r_{i,j} M_{i,j}(i',j')$.
 \item There exist an $n\times n$ bistochastic matrix $R' = (r_{i',j'})$ and $n^2$ bistochastic matrices $M'_{i',j'}$ of dimension $(n-1)\times(n-1)$ whose rows and vertices are indexed by $[n] \setminus \{i'\}, [n] \setminus \{j'\}$ (respectively) such that $M((i,i'),(j,j')) = r_{i',j'} M_{i',j'}(i,j)$.
\end{enumerate}

Every permutation $\pi \in S_n$ gives rise to a 2-bistochastic matrix $M_\pi$ given by $M_\pi((i,i'),(j,j')) = \charf{\pi(i)=j \text{ and } \pi(i') = j'}$. We can now state Theorem~29 in the special case $k=2$.

\begin{thm}{29}
 An $n(n-1) \times n(n-1)$ matrix $M$ is 2-bistochastic if and only if it is a convex combination of $n(n-1) \times n(n-1)$ matrices induced by permutations of $[n]$.
\end{thm}

The following matrix is a counterexample to this theorem, for $n = 4$ (the theorem holds for smaller $n$):
\[
 \begin{array}{c|ccc|ccc|ccc|ccc|}
  &12&13&14&21&23&24&31&32&34&41&42&43\\\hline
12& 0& 1/4& 0& 0& 0& 0& 0& 0& 1/4& 1/4& 1/4& 0\\
13& 0& 0& 1/4& 0& 0& 0& 1/4& 0& 0& 0& 1/4& 1/4\\
14& 1/4& 0& 0& 0& 0& 0& 0& 1/4& 0& 1/4& 0& 1/4\\\hline
21& 0& 0& 1/4& 1/4& 0& 0& 0& 0& 1/4& 0& 0& 1/4\\
23& 0& 1/4& 0& 0& 0& 1/4& 0& 1/4& 0& 1/4& 0& 0\\
24& 1/4& 0& 0& 0& 1/4& 0& 1/4& 0& 0& 0& 1/4& 0\\\hline
31& 0& 0& 1/4& 0& 1/4& 0& 0& 0& 1/4& 1/4& 0& 0\\
32& 1/4& 0& 0& 0& 0& 1/4& 1/4& 0& 0& 0& 0& 1/4\\
34& 0& 1/4& 0& 1/4& 0& 0& 0& 1/4& 0& 0& 1/4& 0\\\hline
41& 0& 1/4& 0& 0& 0& 1/2& 1/4& 0& 0& 0& 0& 0\\
41& 0& 0& 1/4& 1/4& 1/4& 0& 0& 1/4& 0& 0& 0& 0\\
43& 1/4& 0& 0& 1/4& 1/4& 0& 0& 0& 1/4& 0& 0& 0\\\hline
 \end{array}
\]
\[
 \begin{array}{c|ccc|ccc|ccc|ccc|}
  &21&31&41&12&32&42&13&23&43&14&24&34\\\hline
21& 1/4& 0& 0& 0& 0& 0& 0& 0& 1/4& 1/4& 0& 1/4\\
31& 0& 0& 1/4& 0& 0& 0& 0& 1/4& 0& 1/4& 0& 1/4\\
41& 0& 1/4& 0& 0& 0& 0& 1/4& 0& 0& 0& 1/2& 0\\\hline
12& 0& 0& 1/4& 0& 0& 1/4& 1/4& 0& 0& 0& 0& 1/4\\
32& 0& 1/4& 0& 1/4& 0& 0& 0& 0& 1/4& 0& 1/4& 0\\
42& 1/4& 0& 0& 0& 1/4& 0& 0& 1/4& 0& 1/4& 0& 0\\\hline
13& 0& 1/4& 0& 0& 0& 1/4& 0& 0& 1/4& 1/4& 0& 0\\
23& 0& 0& 1/4& 0& 1/4& 0& 1/4& 0& 0& 0& 1/4& 0\\
43& 1/4& 0& 0& 1/4& 0& 0& 0& 1/4& 0& 0& 0& 1/4\\\hline
14& 0& 0& 1/4& 1/4& 1/4& 0& 0& 0& 1/4& 0& 0& 0\\
24& 0& 1/4& 0& 1/4& 0& 1/4& 0& 1/4& 0& 0& 0& 0\\
34& 1/4& 0& 0& 0& 1/4& 1/4& 1/4& 0& 0& 0& 0& 0\\\hline
 \end{array}
\]

The reader can check that this matrix (given in two different orderings of the rows and columns) is 2-bistochastic. According to Theorem~29, it should be a convex combination of matrices induced by permutations. However, its support does not contain the support of any matrix induced by a permutation.

\smallskip

We go on to pinpoint the mistake in the proof of Theorem~29; what follows uses the same notation as the original proof. Let $M$ be a 2-bistochastic matrix. According to condition~(S1), there are bistochastic matrices $R = (r_{i,j})$ and $M_{i,j}$ such that $M((i,i'),(j,j')) = r_{i,j} M_{i,j}(i',j')$. Since $R$ is bistochastic, according to Birkhoff's theorem it is a convex combination of permutation matrices. If $P$ is one of these permutation matrices, then we can write $R = sP + (1-s)T$ for some $s \in (0,1]$. The proof then considers the modified matrix $\tilde{M}$ \emph{(my notation)} defined by $\tilde{M}((i,i'),(j,j')) = P(i,j) M_{i,j}(i',j')$. Under the (tacit) assumption that $\tilde{M}$ is 2-bistochastic, the proof goes on to show that $\tilde{M}$ is the matrix induced by the permutation corresponding to $P$. However, while $\tilde{M}$ certainly satisfies~(S1), it need not satisfy~(S2). For a concrete example, consider the matrix $M$ above. The corresponding matrix $R$ is
\[
\begin{pmatrix}
1/4 &  0  & 1/4 & 1/2 \\
1/4 & 1/4 & 1/4 & 1/4 \\
1/4 & 1/4 & 1/4 & 1/4 \\
1/4 & 1/2 & 1/4 &  0 \\
\end{pmatrix} =
\frac{1}{4}
\begin{pmatrix}
 0 & 0 & 0 & 1 \\
 0 & 0 & 1 & 0 \\
 0 & 1 & 0 & 0 \\
 1 & 0 & 0 & 0
\end{pmatrix} +
\frac{1}{4}
\begin{pmatrix}
 1 & 0 & 0 & 0 \\
 0 & 1 & 0 & 0 \\
 0 & 0 & 0 & 1 \\
 0 & 0 & 1 & 0
\end{pmatrix} +
\frac{1}{4}
\begin{pmatrix}
 0 & 0 & 1 & 0 \\
 0 & 0 & 0 & 1 \\
 1 & 0 & 0 & 0 \\
 0 & 1 & 0 & 0
\end{pmatrix} +
\frac{1}{4}
\begin{pmatrix}
 0 & 0 & 0 & 1 \\
 1 & 0 & 0 & 0 \\
 0 & 0 & 1 & 0 \\
 0 & 1 & 0 & 0
\end{pmatrix}.
\]
Denote the first permutation matrix by $P$. If we replace $R$ by $P$ then we get the following matrix:
\[
 \begin{array}{c|ccc|ccc|ccc|ccc|}
  & 12 & 13 & 14 & 21 & 23 & 24 & 31 & 32 & 34 & 41 & 42 & 43 \\\hline
 12 & 0 & 0 & 0 & 0 & 0 & 0 & 0 & 0 & 0 & 1/2 & 1/2 & 0 \\
 13 & 0 & 0 & 0 & 0 & 0 & 0 & 0 & 0 & 0 & 0 & 1/2 & 1/2 \\
 14 & 0 & 0 & 0 & 0 & 0 & 0 & 0 & 0 & 0 & 1/2 & 0 & 1/2 \\\hline
 21 & 0 & 0 & 0 & 0 & 0 & 0 & 0 & 0 & 1 & 0 & 0 & 0 \\
 23 & 0 & 0 & 0 & 0 & 0 & 0 & 0 & 1 & 0 & 0 & 0 & 0 \\
 24 & 0 & 0 & 0 & 0 & 0 & 0 & 1 & 0 & 0 & 0 & 0 & 0 \\\hline
 31 & 0 & 0 & 0 & 0 & 0 & 0 & 0 & 1 & 0 & 0 & 0 & 0 \\
 32 & 0 & 0 & 0 & 0 & 0 & 0 & 0 & 0 & 1 & 0 & 0 & 0 \\
 34 & 0 & 0 & 0 & 1 & 0 & 0 & 0 & 0 & 0 & 0 & 0 & 0 \\\hline
 41 & 0 & 1 & 0 & 0 & 0 & 0 & 0 & 0 & 0 & 0 & 0 & 0 \\
 42 & 0 & 0 & 1 & 0 & 0 & 0 & 0 & 0 & 0 & 0 & 0 & 0 \\
 43 & 1 & 0 & 0 & 0 & 0 & 0 & 0 & 0 & 0 & 0 & 0 & 0 \\\hline
 \end{array}
\]
\[
 \begin{array}{c|ccc|ccc|ccc|ccc|}
  & 21 & 31 & 41 & 12 & 32 & 42 & 13 & 23 & 43 & 14 & 24 & 34 \\\hline
 21 & 0 & 0 & 0 & 0 & 0 & 0 & 0 & 0 & 0 & 0 & 0 & 1 \\
 31 & 0 & 0 & 0 & 0 & 0 & 0 & 0 & 1 & 0 & 0 & 0 & 0 \\
 41 & 0 & 0 & 0 & 0 & 0 & 0 & 1 & 0 & 0 & 0 & 0 & 0 \\\hline
 12 & 0 & 0 & 1/2 & 0 & 0 & 1/2 & 0 & 0 & 0 & 0 & 0 & 0 \\
 32 & 0 & 0 & 0 & 0 & 0 & 0 & 0 & 0 & 0 & 0 & 1 & 0 \\
 42 & 0 & 0 & 0 & 0 & 0 & 0 & 0 & 0 & 0 & 1 & 0 & 0 \\\hline
 13 & 0 & 0 & 0 & 0 & 0 & 1/2 & 0 & 0 & 1/2 & 0 & 0 & 0 \\
 23 & 0 & 0 & 0 & 0 & 1 & 0 & 0 & 0 & 0 & 0 & 0 & 0 \\
 43 & 0 & 0 & 0 & 1 & 0 & 0 & 0 & 0 & 0 & 0 & 0 & 0 \\\hline
 14 & 0 & 0 & 1/2 & 0 & 0 & 0 & 0 & 0 & 1/2 & 0 & 0 & 0 \\
 24 & 0 & 1 & 0 & 0 & 0 & 0 & 0 & 0 & 0 & 0 & 0 & 0 \\
 34 & 1 & 0 & 0 & 0 & 0 & 0 & 0 & 0 & 0 & 0 & 0 & 0 \\\hline
 \end{array}
\]
As can be seen, this matrix satisfies (S1) but not (S2).

\pagebreak

\section{Final remarks} \label{sec:remarks}

The counterexample described to Theorem~27 arose while trying to prove a similar theorem on the \emph{slice} $J(n,k) = \{ x \in \{0,1\}^n : \sum_i x_i = k \}$. Specifically, we were interested in understanding the structure of Boolean functions of degree~$d$ (the degree of a function on $J(n,k)$ is the minimal degree of a polynomial defining it). Every Boolean function $f$ on the slice $J(n,k)$ can be lifted to a Boolean function $F$ on $S_n$ using the formula
\[
 F(\pi) = f(\pi(\{1,\ldots,k\})).
\]
Moreover, if $f$ has degree $d$ then $F \in V_d$. It is then a simple exercise to deduce the following theorem from Theorem~27:

\begin{thm}{27J}
 If $f$ is a nonnegative function on $J(n,k)$ of degree $d$ then $f$ can be written as a nonnegative combination of degree-$d$ monomials in the variables $x_1,\ldots,x_n,1-x_1,\ldots,1-x_n$. Moreover, if $f$ is a Boolean function on $J(n,k)$ of degree $d$ then $f$ can be written as a sum of degree-$d$ monomials (in the same variables), which are necessarily ``disjoint'' in the sense that no two can evluate to~1 at the same time.
\end{thm}

This theorem fails for the following function, which is where our counterexample to Theorem~27 comes from:
\[
 f = \charf{x_1 = x_2 = x_3} = 1 - x_1 - x_2 - x_3 + x_1x_2 + x_1x_3 + x_2x_3.
\]
Curiously, it does hold for the negation of $f$:
\[
 1-f = x_1(1-x_2) + x_2(1-x_3) + x_3(1-x_1).
\]
As noticed by David Ellis, this example can be extended to any even degree.

It is natural to ask whether Theorem~27 or Theorem~27J can be corrected. Here are two possible ways to correct Theorem~27J, in the special case of Boolean functions:
\begin{enumerate}
 \item For every $d$ there is a constant $s_d$ such that every Boolean function on $J(n,k)$ of degree $d$ can be written as a sum monomials of degree~$s_d$.
 \item For every $d$ there is a constant $r_d$ such that for every Boolean function $f$ on $J(n,k)$ of degree $d$, either $f$ or $1-f$ can be written as a sum of monomials of degree~$r_d$.
\end{enumerate}

For a Boolean function $f$ on $\{0,1\}^n$, the minimal $s$ such that $f$ can be written as a sum of monomials of degree~$s$ is known as its one-sided \emph{unambiguous certificate complexity}, which is known to be polynomially related to the degree~\cite{Goos,BenDavid}. The analog of the parameter $r$ is known as the two-sided unambiguous certificate complexity, and is also polynomial related to the degree~\cite{Goos,BenDavid}. Therefore for the Boolean cube we know that $s_d,r_d = d^{\Theta(1)}$.

A classical example showing that $r_2 > 2$ (on both the Boolean cube and $J(n,k)$ for appropriate $n,k$) is the 4-sortedness function
\[ f(x_1,x_2,x_3,x_4) = \charf{x_1 \leq x_2 \leq x_3 \leq x_4 \text{ or } x_1 \geq x_2 \geq x_3 \geq x_4}. \]
While not apparent from this definition, $f$ has degree~2. It is not hard to check that the unambiguous certificate complexity of both $f$ and $1-f$ is larger than~2. By iterating this function, we get the lower bound $r_d = d^{\Omega(1)}$.

Finally, let us mention, that an approximate version of Theorem~27 does hold for sparse Boolean functions, namely the main result of~\cite{EFF3}. Roughly speaking, this result states that if $f$ is a Boolean function close to $V_d$ (in $L^2$) and $\mathbb{E}[f] = O(1/n^d)$ then $f$ is close to a union of $d$-cosets. For similar results in the case $d=1$, see~\cite{EFF1,EFF2}.
	
\bibliographystyle{alpha}
\bibliography{EFP-comment}

\end{document}